\RequirePackage{fix-cm}
\documentclass[smallextended]{svjour3}       
\smartqed  
\usepackage{graphicx}
\usepackage{color}
\journalname{Mathematical Intelligencer}
\begin{document}

\title{Dynamical systems, celestial mechanics, and music: \\ Pythagoras revisited}

\author{Julyan H. E. Cartwright    \and
          Diego L. Gonz\'alez  \and 
          Oreste Piro
}


\institute{
Julyan H. E. Cartwright \at  Instituto Andaluz de Ciencias de la Tierra CSIC--Universidad de Granada \\ Armilla, 18100 Granada \\ Spain \\ and \\ Instituto Carlos I de F\'{\i}sica Te\'orica y Computacional Universidad de Granada \\ 18071 Granada, Spain  \\
\email{julyan.cartwright@csic.es}
\and
Diego L. Gonz\'alez \at
              Istituto per la Microelettronica e i Microsistemi Area della Ricerca CNR di Bologna \\ 40129 Bologna, Italy  \\
                 and \\ Dipartimento di Scienze Statistiche ``Paolo Fortunati'' Universit\`a di Bologna \\ 40126 Bologna, Italy  \\
              \email{gonzalez@bo.imm.cnr.it}           
           \and
           Oreste Piro \at
  Departament de F\'{\i}sica, Universitat de les Illes Balears, \\ 07071 Palma de Mallorca, Spain \\
  \email{oreste.piro@uib.es}
  }
  
\date{Received: date / Accepted: date}

\maketitle

\begin{quote}
``I am every day more and more convinced of the Truth of Pythagoras's Saying, that Nature is sure to act consistently, and with a constant Analogy in all her Operations: From whence I conclude that the same Numbers, by means of which the Agreement of Sounds affects our Ears with Delight, are the very same which please our Eyes and Mind. We shall therefore borrow all our Rules for the Finishing our Proportions, from the Musicians, who are the greatest Masters of this Sort of Numbers, and from those Things wherein Nature shows herself most excellent and compleat.''
\end{quote}
[Leon Battista Alberti, (1407--1472)  De Re Aedificatoria. Chapter V of Book IX of Alberti's Ten Books of Architecture (James Leoni, translator).]

\newpage 

\noindent
Gioseffo Zarlino  reintroduced the Pythagorean paradigm into Renaissance musical theory. In a similar fashion,  Nicolaus Copernicus, Galileo Galilei, Johannes Kepler, and Isaac Newton reinvigorated Pythagorean ideas in celestial mechanics; Kepler and Newton explicitly invoked musical principles. Today, the theory of dynamical systems allows us to describe very different applications of physics, from the orbits of asteroids in the Solar System to the pitch of complex sounds. 
Our aim in this text is to review the overarching aims of our research in this field over the past quarter of a century.
We demonstrate with a combination of dynamical systems theory and music theory the thread running from Pythagoras to Zarlino that allowed the latter to construct musical scales using the ideas of proportion known to the former, and we discuss how the modern theory of dynamical systems, with the study of resonances in nonlinear systems, returns  to Pythagorean ideas of a Musica Universalis.


\section*{Nature most excellent and compleat}

To the Pythagoreans, music represented the paradigm of order emerging from the primordial chaos, that is, the root of the cosmos. Moreover, this order could be understood through number, the tool used by God, the Great Geometer, for creating the universe. For, if music were governed by number --- Pythagorean proportion --- the same was true of the celestial bodies that play a heavenly symphony in the sky: the music of the spheres. At the same time, on an intermediate level, the temple represented the \emph{trait d'union} between man and the heavens, the micro- and macro-cosmos. Within this metaphysical framework it is not at all strange that, as discussed by Alberti \cite{alberti}, the architecture of the temple should also follow the laws of numbers, that is, the particular proportions giving a sense of harmony and beauty through all the levels of creation. 

The Pythagorean approach lay almost dormant during mediaeval times. However, beginning in the High Middle Ages, a  revival began to take shape. At the beginning of the fourteenth century, Dante \cite[II-xiii-18]{dante} wrote that regarding the three principles of natural things, namely matter, privation, and form:
\begin{quote}
Number exists not only in all of them together, but also, upon careful reflection, in each one individually; for this reason Pythagoras, as Aristotle says in the first book of the Physics, laid down even and odd as the principles of natural things, considering all things to have numerical aspect.
\end{quote}
With the advent of the Renaissance, Pythagorean ideas spread  with new vigour.  Many of the translations of Arabic books that had preserved classical Greek knowledge through mediaeval times  took place in Venice where, thanks to maritime commerce, translators from the Arabic and  from the Greek were available for producing Latin and vulgar versions of the texts. It is natural, then, that the reintroduction of Pythagorean ideas in musical theory and practice found full expression in Venice. The main protagonist was Gioseffo Zarlino, who consolidated the earlier work of Franchino Gaffurio and Francesco Maurolico.
Zarlino, a Franciscan friar, was organist in Chioggia and then choir-master at St. Mark's in Venice from 1565 to 1590. 
 In architecture, its main standard-bearer was Leon Battista Alberti, quoted above, whose work was further developed by many of the great eclectic minds of the Renaissance, including Francesco de Giorgio Martini, Sandro Botticelli, Leonardo da Vinci, Francesco Giorgi and Andrea Palladio.
 
With the nascence of modern experimental science following the work of Galileo Galilei --- intellectual descendent of Nicolaus Copernicus, contemporary of Johannes Kepler, and academic forefather of Isaac Newton --- there began a gradual process of bifurcation between arts and sciences. 
On one hand, there has been a progressive demystification of the role played by proportion in the arts, including in music, architecture, 
and painting. 
On the other hand, this has been 
accompanied by an increasing neglect of science for the Pythagorean ideas relating art,  music, and natural philosophy. 

However, recent research in dynamical systems theory should change these views. Dynamical resonances predicted by  theory describe very different types of behaviour, from the pitch of complex sounds to the orbits of celestial bodies in the Solar System. These achievements of modern science bring together in a surprising fashion physiological behaviour, astronomy, and the theory of numbers; that is, the micro-cosmos, the macro-cosmos, and Pythagorean number. 
In this text we discuss how the mathematics of proportion and aesthetics that developed from Pythagoras to the Renaissance has re-emerged in the science of dynamical systems. We highlight how modern theory allows new insight into the Pythagorean scientific conception reintroduced into music by Zarlino and others in Renaissance Venice.

\section*{Pythagoras, musical intervals, and the Music of the Spheres}

\begin{figure}
  \begin{center}
 \includegraphics*[width= 0.49\textwidth]{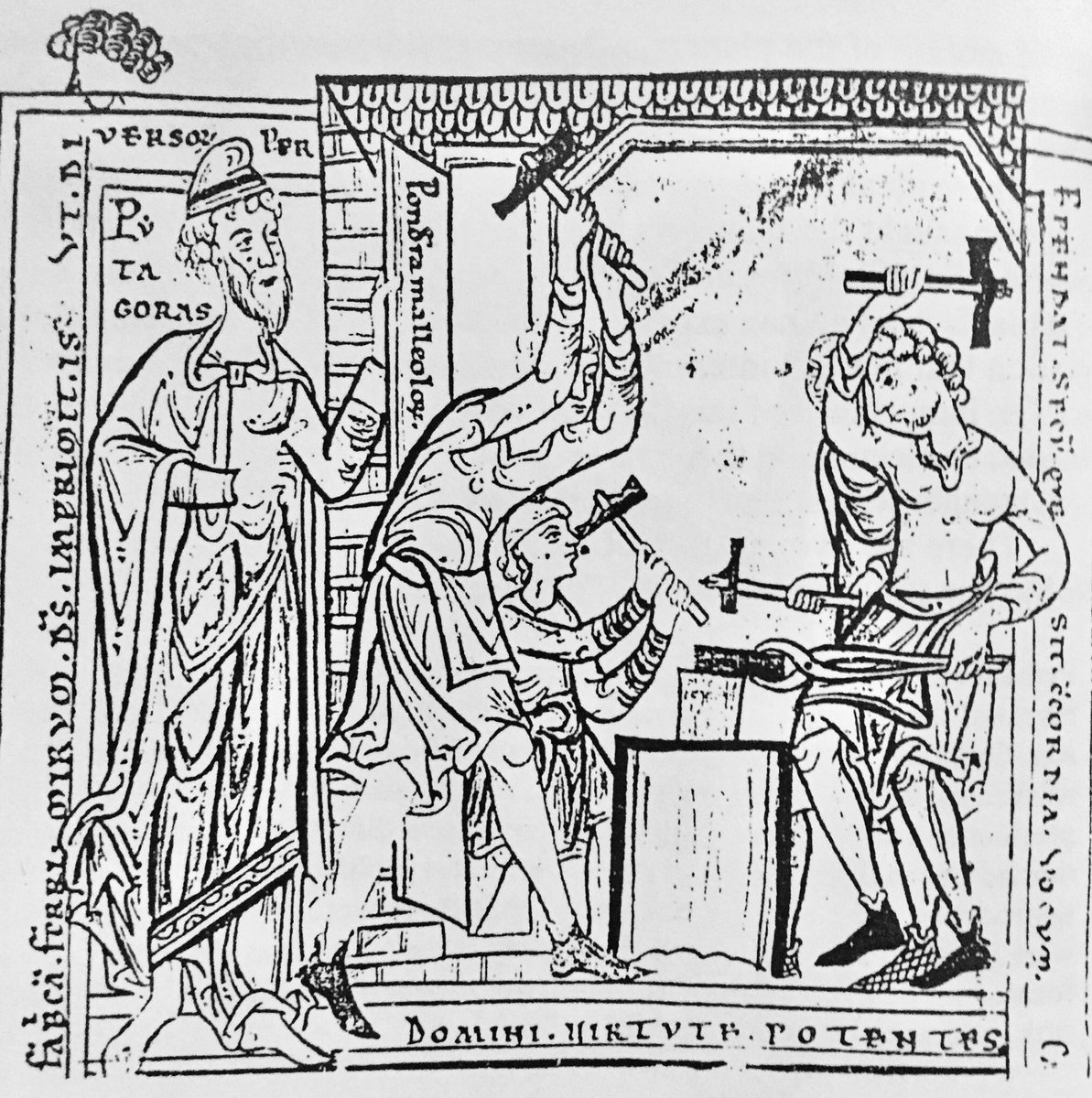}\\
  \caption{Pythagoras and the discovery of musical intervals in a blacksmith's forge illustrated in a mediaeval woodcut (Bayerische Staatsbibliothek). Today, Verdi's \emph{Il Trovatore} and Wagner's \emph{Das Rheingold} and \emph{Siegfried}, among other works, all use the sounds of hammers on anvils. But real hammers on real anvils are notoriously nonmusical, so productions of these operas often use fabricated percussion instruments made to look anvil-like that do emit musical tones.  As was already pointed out by scholars of musical history some centuries ago ``upon examination and experiment it appears, that hammers of different size and weight will no more produce different tones upon the same anvil, than bows or clappers of different sizes will from the same string or bell. Indeed, both the hammers and anvils of antiquity must have been of a construction very different from those of our degenerate days, if they produced any tones  that were strictly musical'' \cite{burney1776}.
 Jones  \cite{jones1818} (page 344) makes the  suggestion that ancient Greeks might have used convex shield-like pieces of metal, which might possibly have been ``sonorous anvils", and goes on to propose how this might resolve the question, but we are not convinced.
  }
  \label{fig:pythagoras}
  \end{center}
\end{figure}

For the Pythagorean school \cite{riedweg2005,zhmud2012}, the Cosmos, that is, our universe, that is, was nothing but the result of the order imposed by the Demiurge, the Great Geometer, on the primitive chaos. The identification of the god of creation with a great surveyor shows clearly the concepts that guided Pythagorean thought. 
The Pythagoreans saw that there is a unique tool to find order in the universe: mathematics, which may occur in its twin aspects, of geometry and arithmetic. 
These two, together with music and astronomy, became the four liberal arts of the quadrivium. The quadrivium brought about the gradual introduction of classical thought into mediaeval instruction. Its consolidation is attributed to Boethius in the 6th century,   a main mediaeval intermediary from the Pythagoreans and classical antiquity to the neo-Platonists of the Renaissance. His near contemporary Proclus wrote \cite{proclus}
\begin{quote}
	The Pythagoreans considered all mathematical science to be divided into four parts: one half they marked off as concerned with quantity, the other half with magnitude; and each of these they posited as twofold. A quantity can be considered in regard to its character by itself or in its relation to another quantity, magnitudes as either stationary or in motion. Arithmetic, then, studies quantities as such, music the relations between quantities, geometry magnitude at rest, spherics [i.e., astronomy] magnitude inherently moving.
\end{quote}
The properties of numbers were the most important subject of study since these properties could be observed at all levels of Creation, from the movement of the stars in the firmament --- the macrocosm --- down to man himself --- the microcosm. It is within this context that the scientific study of music began. In the famous legend of the forge \cite{anderson1983}, first related by Nicomachus \cite{nicomachus}, illustrated in Figure~\ref{fig:pythagoras}, the discovery of harmonic musical intervals is attributed to Pythagoras himself as he listened to the sounds produced by hammers of different sizes that struck a large piece of red-hot iron on an anvil. Although this tradition is based on ancient myths and legends --- which notably associate the musical sound of blacksmiths at work and the invention of two sciences: acoustics and metallurgy --- most likely the Pythagoreans studied musical intervals not with hammers, but with the monochord and other stringed instruments, as we note in Fig.~\ref{fig:pythagoras}. Pythagoras' discovery is probably the first in human history that can be qualified as a scientific theory; i.e., a description of a natural phenomenon in mathematical terms. With the finding that certain small integer relationships between the lengths of  strings (see Table~\ref{table:intervals}) produce harmonious sounds, the Pythagoreans put music at the centre of the intellectual effort of their school.
In particular this was due to Archytas, ``of all the Pythagoreans the most devoted to the study of music'', who ``tried to preserve what follows the principles of reason not only in the concords but also in the divisions of the tetrachords'', according to Ptolomy \cite{ptolemy,barker1994}.
This view of music  as a paradigm of cosmic order crystallized completely with the music of the spheres: the idea that the planets through their movement in the sky produce subtle harmonies that only the ears of initiates --- such as Pythagoras --- could hear.

\begin{table}
\caption{
The consonant intervals ---  ratios of frequencies between  two tones --- of modern Western music and the principal intervals of the so-called just scales are characterized by rational numbers. 
For example,
when the frequency of a C is multiplied by 3, one gets a (just) G one octave higher. To get the interval C--G within an octave, one must divide the larger one by 2, that is 3/2. Likewise dividing the frequency of C by 3, one obtains the (just) F two octaves lower, and one must multiply that frequency by 4 to obtain the F in the octave above the C, giving the interval 4/3. 
}\label{table:intervals}
\centering\begin{tabular}{ccc}
\hline\hline
Intervals &	Quality &	Rational number  \\
\hline
Unison &	Consonant &	1/1 \\
Octave &&	2/1 \\
Major Sixth &&	5/3 \\
Minor Sixth &&	8/5 \\
Fifth &&	3/2 \\
Fourth &&	4/3 \\
Major Third &&	5/4 \\
Minor Third &&	6/5 \\
Major Tone &	Dissonant	& 9/8 \\
Minor Tone &&	10/9 \\
Major Semitone &&	16/15 \\
Minor Semitone &&	25/24 \\
Syntonic Comma &&		81/80 \\
\hline\hline
\end{tabular}
\end{table}

Pythagorean thought in its esoteric key was greatly demystified by the birth of modern experimental science with Galileo Galilei. However, the concepts developed  by the adepts of the Pythagorean school clearly represent an encrypted or at least an analogous version of the reality they were trying to describe.

\section*{Zarlino and the reintroduction of Pythagorean concepts in the music of Renaissance Venice}

Gioseffo Zarlino is recognized as the main protagonist in the reintroduction of classical concepts into the theory of Western music. 
Zarlino produced important theoretical papers that were compiled in several volumes: \emph{Le Istitutioni Harmoniche} of 1558, \emph{Dimonstrationi Harmoniche} of 1571, and \emph{Sopplimenti Musicali} of 1588, all published in Venice; they were collected in his Complete Works \cite{Zarlino1589}. We must note here the disputes in the field of music between Zarlino and his pupil, Vincenzo Galilei, musician and the father of Galileo Galilei \cite{galilei1589,drake1970,walker1973}, which sound like  a current-day bad-tempered argument between a theoretician and an experimentalist. 

\begin{figure}
  \begin{center}
 \includegraphics*[width=0.8\textwidth]{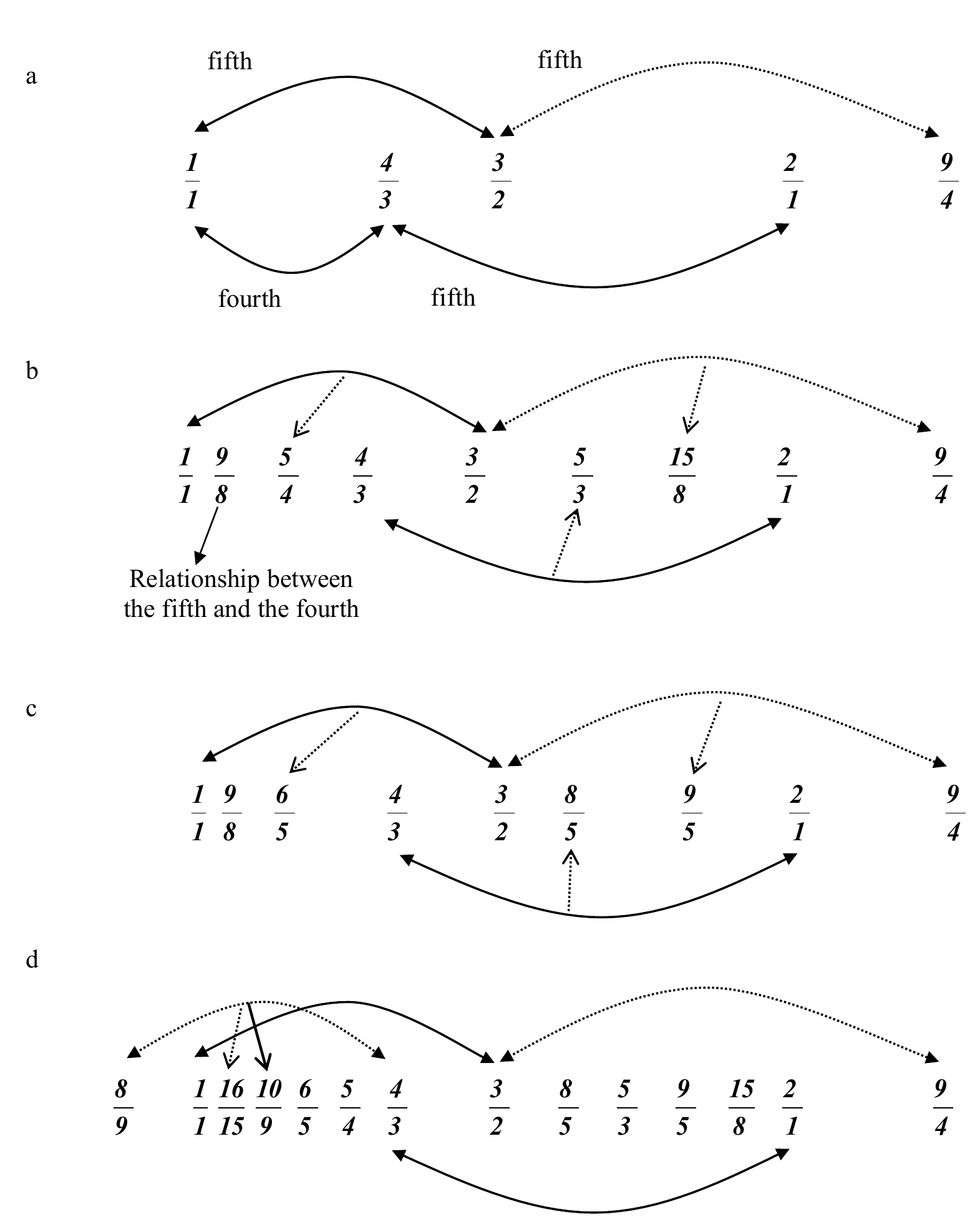}\\
  \caption{a) Pythagorean intervals: the fifth and fourth can be obtained by calculating, respectively, the arithmetic and harmonic means of the octave (double arrows indicate the interval between notes (fifth or fourth), dotted double arrows indicate that the corresponding interval is formed with a note outside the octave); b) Zarlino major scale: adding to the notes represented in the previous figure the arithmetic means of the fifths represented in the same figure, and the whole tone $9/8$ which is obtained as the ratio between the fifth and fourth, one obtains the Zarlino major scale (dotted arrows indicate arithmetic mean); c) Zarlino minor scale: although there are several variants to obtain a chromatic scale starting from basic intervals (for example the variant of Delezenne or the use of a minor semitone $25/24$ to obtain the altered notes (\cite{Lattard1988}, p. 24)), a Zarlino minor scale can be obtained in a single step by calculating the harmonic means of the fifths shown in a) (dotted arrows indicate harmonic mean); d) Convergence with the golden scale of 12 notes: As mentioned in the text, to get all the possible fifths starting from the Pythagorean fifths and fourths in the octave is necessary to include the note $8/9$ (a fifth below $4/3$). If one calculates the arithmetic and harmonic means of this fifth interval one obtains the notes $16/15$ and $10/9$.
 If they are put into place, we obtain all the rational intervals of the 12-note golden scale (a dotted arrow indicates harmonic mean, a plain arrow indicates arithmetic mean)  \cite{Cartwright2002,Gonzalez2002}.
}
  \label{fig:intervals}
  \end{center}
\end{figure}

Zarlino widened the range of Pythagorean harmonic intervals so as to define the major and minor scales that are today known as \emph{just};  that is, scales in which all the intervals are defined by rational numbers. The inclusion of the third and sixth intervals, both major and minor, defined by rational relationships like the Pythagoreans' fourth and fifth but extended to the senarius --- the group of numbers 1--6 --- permitted the generation of the major and minor scales.
His Sintono Diatonico, which is not his original invention, but in fact dates back to Ptolemy in the second century CE, represents a major step towards modern polyphonic music in anticipating the theory of harmonics of Sauveur \cite{sauveur} of the eighteenth century. The influence of the work of Zarlino spread throughout Europe in the succeeding centuries, and led to the diffusion of concepts like just scales, developed rigorously from  Helmholtz \cite{helmholtz1863} in the nineteenth century to current microtonal music. 

It is a mistake to regard  microtonality as modern, however. 
Microtonality --- that is, music implemented with musical scales using intervals smaller than a semitone --- has historically characterized  the music of a number of Asian cultures, but also in the European tradition  there have been many examples.
Already in \emph{Le Istitutioni Harmoniche} Zarlino had represented a harpsichord with 19 notes per octave \cite{Wraight2003}. This proposal of Zarlino, and others, such as 24 or 31 notes per octave, highlights the convergence with microtonal sub-divisions suggested by recent theoretical results that we shall discuss below.
From the theoretical point of view, Zarlino extended the set of the harmonic intervals using the same approach developed by the Pythagorean school:  that of proportional division. The intervals found by Pythagoras are those of the octave, fifth, fourth, and whole tone, which are summarized in Table~\ref{table:intervals}. All these intervals can be obtained by division of the octave using the arithmetic and  harmonic means. 

The ancient Greeks knew of the division of an interval into three different types of means between two quantities, $a$ and $b$:
\begin{enumerate}
\item The arithmetic mean, or simply the mean; defined as 
$$
\langle a, b \rangle = \frac{a+b}{2},
$$
\item the harmonic mean, 
$$
| a, b | = \frac{1}{\displaystyle\left\langle \frac{1}{a}, \frac{1}{b} \right\rangle} = \frac{2ab}{a+b},
$$
equal to the inverse of the arithmetic mean of the inverses,
\item and the geometric mean
$$
\| a, b \| = \sqrt{a b},
$$
which coincides with the square-root of the product of the arithmetic mean and the harmonic mean.
\end{enumerate}
In other words, the geometric mean of the interval is  the geometric mean of the arithmetic and harmonic means of the same interval, a result attributed to Archytas.
Applying the arithmetic mean and the harmonic mean to the octave interval $[1/1, 2/1]$ one obtains the values $3/2$ and $4/3$ corresponding to the fifth and the fourth, while the ratio between these two produces the interval corresponding to the whole tone, $9/8$. The geometric mean was not used for the calculation of harmonic intervals, since its application usually produces irrational numbers, and so it was `imperfect' from the perspective of Pythagorean doctrine.

From a theoretical point of view,  Zarlino's intervals can be obtained by calculating the arithmetic and harmonic means of fifths appearing in Pythagorean subdivisions.
This construction, together with the scale it generates, is represented in Figure~\ref{fig:intervals}.
The construction begins with the recognized Pythagorean intervals that are the octave and fifth.
We may start with an octave of C and assume that in addition we know only of fifths. In other words, we have a C and another C an octave above (1/1; 2/1). Add a fifth to the first C. It is the note G (3/2), but now we can insert a fifth under the upper C, which is the note F (3/4). This note in relation to the lower C gives a fourth. In this way we generate F and G that can be defined as intervals of fifths with the two Cs forming  the octave.
As we use only fifths, the note G can generate a fifth  above, which is 9/4, and the F a fifth below, which is 8/9.
Let us now remain with the initial notes C to C (1/1, 2/1), and their upper and lower fifths, F, G, D and B$\flat$; 4/3, 3/2, 9/4, and 8/9.
While the initial Pythagorean scheme uses only the numbers from 1 to 4, or their powers, Zarlino extends this to include 5 and 6; in this way one can include the intervals of sixths and thirds.
These numbers appear naturally as different means of the existing intervals, for example --- Figure~\ref{fig:intervals}b --- including the arithmetic means gives the major scale of Zarlino (obtained without performing the arithmetic mean with the lower fifth of the F (8/9); it is necessary to include as a note the naturally generated interval between the Pythagorean F and G, this is a whole tone, 3/2 / 4/3 = 9/8).
It is enough  now to include the harmonic means to obtain the Zarlino minor scale (again the fifth lower than the F is not necessary); Figure~\ref{fig:intervals}c.
Thus, starting only with the notes of the Pythagorean intervals, the octave and fifth, and the fourth and the whole tone as a consequence, and adding the arithmetic and harmonic means, the major scale and the minor scale of Zarlino are obtained.

\section*{Science's return to Pythagoras through celestial mechanics}

The Pythagorean \emph{Musica Universalis} or music of the spheres is today a metaphor, but the historical connections between music, mathematics, and astronomy have had a profound impact upon all these disciplines \cite{archibald1924,proust2011}. 
In the preface to his book \emph{De revolutionibus orbium coelestium} (On the Revolutions of the Heavenly Spheres, \cite{copernicus1543}), Nicolaus Copernicus cites Pythagoreans as the most important influences on the development of his heliocentric model of the universe. In the century following his death the Inquisition banned \emph{De revolutionibus} and all books advocating the Copernican system, which it called ``the false Pythagorean doctrine, altogether contrary to Holy Scripture.'' This ban was flouted by Galileo Galilei in his \emph{Dialogo sopra i due massimi sistemi del mondo}  (Dialogue Concerning the Two Chief World Systems, \cite{galileo1632}) with its account of conversations between a Copernican scientist, Salviati, who represents the author, a witty scholar, Sagredo, and a plodding Aristotelian, Simplicio. Galileo has Salviati say 
\begin{quote}
That the Pythagoreans held the science of numbers in high esteem, and that Plato himself admired the human understanding and believed it to partake of divinity simply because it understood the nature of numbers, I know very well; nor am I far from being of the same opinion. But that these mysteries which caused Pythagoras and his sect to have such veneration for the science of numbers are the follies that abound in the sayings and writings of the vulgar, I do not believe at all. Rather I know that, in order to prevent the things they admired from being exposed to the slander and scorn of the common people, the Pythagoreans condemned as sacrilegious the publication of the most hidden properties of numbers or of the incommensurable and irrational quantities which they investigated. They taught that anyone who had revealed them was tormented in the other world. Therefore I believe that some one of them, just to satisfy the common sort and free himself from their inquisitiveness, gave it out that the mysteries of numbers were those trifles which later spread among the vulgar.
\end{quote}
The work begins
\begin{quote}
Several years ago there was published in Rome a salutary edict which, in order to obviate the dangerous tendencies of our present age, imposed a seasonable silence upon the Pythagorean opinion that the Earth moves. There were those who impudently asserted that this decree had its origin not in judicious inquiry, but in passion none too well informed. Complaints were to be heard that advisers who were totally unskilled at astronomical observations ought not to clip the wings of reflective intellects by means of rash prohibitions.
\end{quote}
It was when the inevitable occurred and Galileo was convicted by the Inquisition that he is famously held to have muttered
``eppur si muove'' (``yet it does move'').

\begin{figure}
  \begin{center}
 \includegraphics*[width= 0.6\textwidth]{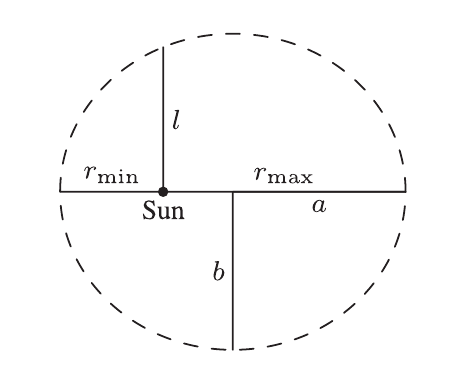}\\
  \caption{Means in a Keplerian elliptical orbit: $a$ is the arithmetic mean, $b$ is the geometric mean, and $l$ is the harmonic mean of $r_{\rm min}$ and $r_{\rm max}$ \cite{mahajan2019}.
}
  \label{fig:kepler_orbit}
  \end{center}
\end{figure}

His contemporary Johannes Kepler's search for harmonic proportions in the Solar System led to his discovery of the laws of planetary motion. Kepler's initial studies in the subject consisted in attempting to reconcile planetary orbits with the geometries of the five Platonic solids; thus he placed an octahedron between the orbits of Mercury and Venus, an icosahedron between Venus and Earth, a dodecahedron between Earth and Mars, a tetrahedron between Mars and Jupiter and a cube between Jupiter and Saturn. This first result he published in \emph{Mysterium Cosmographicum} (Cosmographic Mystery, \cite{kepler1596}). 
He entitled his later book on the subject \emph{Harmonices Mundi} (Harmonics of the World, \cite{kepler1619}), after the Pythagorean teaching that had inspired him. Kepler describes himself falling asleep to the sound of the heavenly music, ``warmed by having drunk a generous draught ... from the cup of Pythagoras''.
Kepler asked whether the greatest and least distances between a planet and the Sun (aphelion and perihelion) might approximate any of the harmonic ratios, but found they did not.
He then looked at the relative maximum and minimum angular velocities of the planets (at perihelion and aphelion) measured from the Sun, and found that planets did seem to approximate harmonic proportions with respect to their own orbits, allowing them to be allotted musical intervals.
He found that the maximum and minimum speeds of Saturn differed by an almost perfect 5/4 ratio, a major third.  The extreme motions of Jupiter differed by 6/5, a minor third. The extremal speeds of Mars, Earth, and Venus approximated 3/2, a fifth; 16/15, a major semitone; and 25/24, a minor semitone,  respectively.
Kepler further examined the ratios between the fastest and slowest speeds of a planet and those of its neighbours, which he called converging and diverging motions. 
His finding that harmonic relationships structure the characteristics of the planetary orbits individually, and their relationships to one another, led to his three laws of planetary motion.
It is beautiful, as we have seen pointed out very recently \cite{mahajan2019}, that Keplerian planetary orbits display all three means discussed above (Fig.~\ref{fig:kepler_orbit}). The semimajor axis $a$ is the arithmetic mean of the perihelion  $r_{\rm min}$ and the aphelion  $r_{\rm max}$,  their geometric mean is the semimajor axis $b$, and $l$, the semilatus rectum, is the harmonic mean.
It is the semimajor axis $a$ whose cube is proportional to the square of the orbital period in Kepler's third law.

Isaac Newton, the discoverer of the secrets of gravity, characteristically did not acknowledge his debt to his immediate predecessors 
but saw in the Pythagorean music of the spheres a prior description of his own law of universal gravitation.  A visitor wrote that
``Mr.\ Newton believes that he has discovered quite clearly that the ancients like Pythagoras, Plato etc.\ had all the demonstrations that he has given on the true System of the World'' \cite{McGuire1966}.
To give the same note (pitch) in the vibration of two strings of different lengths and held by different masses but otherwise similar, according to the law discovered by Mersenne \cite{Mersenne1637}  it is necessary to vary the masses hanging from them in direct relation to square of their lengths. If we now think of two bodies at different distances from a centre of gravitational attraction, like two planets revolving around the Sun, to achieve the same force of attraction it is necessary for the masses to vary following this same law. That is, they should grow as the square of the distance to the centre of attraction. The subject has been discussed by several authors (see, for example, \cite{Gouk1988}, p. 120), but invoking erroneously a dependency with the inverse square of the distance, which defines, instead, the variation of the gravitational force according to Newton's laws.  It was Newton's pupil MacLaurin \cite{Maclaurin1775} who expressed the matter clearly: 
\begin{quote}
In general, that any musical chord may become unison to a lesser chord of the same kind, its tension must be increased in the same proportion as the square of its length is greater; and that the gravity of a planet may become equal to the gravity of another planet nearer to the sun, it must be increased in proportion as the square of its distance from the sun is greater. If therefore we should suppose musical chords extended from the sun to each planet, that all these chords might become unison, it would be requisite to increase or diminish their tensions in the same proportions as would be sufficient to render the gravities of the planets equal. And from the similitude of those proportions the celebrated doctrine of the harmony of the spheres is supposed to have been derived. 
\end{quote}
Probably it is because this was enunciated clearly in print by the disciple, not by Newton himself,  that little attention has been paid to this equivalence between gravitation and music.

Recent results of scientific research in fields from particle physics to cosmology to nonlinear dynamics  return us in curious and interesting ways to Pythagorean ideas. String theory, postulated as a unified description of gravity and particle physics, in which matter is hypothesized to consist at its lowest level of immensely tiny vibrating strings, takes us straight back to Pythagoras. The big-bang theory for the origin of the universe does not differ greatly from the Pythagorean myth of creation: as for the Pythagoreans, for modern physics too there exists something before the Cosmos, the so-called vacuum fluctuations, which are the equivalent of the Pythagorean primordial Chaos, and it is this chaos that is ordered, giving rise to the creation of matter in the universe. In addition, the distribution of the background radiation in the universe, which reflects the structure of the early universe, seems to be due to the propagation of a primordial sound in a manner Pythagoras might well have appreciated \cite{Hu2004}, and a Platonic geometry of the universe has even been suggested \cite{Luminet2003}. 

As we shall describe in more detail below, both celestial mechanics and music theory have a mathematical basis in terms of dynamical systems. 
The existence of a special relationship defined by whole numbers and harmonic intervals as small as those of music can be associated with the properties of nonlinear dynamical systems.
In particular, harmonic musical intervals seen as nonlinear responses of our sense of hearing \cite{Cartwright1999a,Cartwright2001,Gonzalez2002}, and the synchronized movements of the celestial bodies of our planetary system \cite{Escribano2008,Vanyo2011} may be described by resonances of two or more frequencies \cite{Cartwright1999b,Calvo2000,Cartwright2010} in much the same way as first suggested by Pythagoras some 2500 years ago.

\section*{The optimal number of notes in a scale}

\begin{figure}
  \begin{center}
 \includegraphics*[width= 0.6\textwidth]{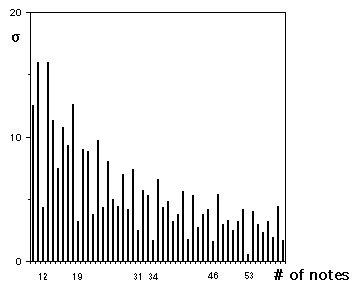}\\
  \caption{Measure of the quality of the approximation of all Zarlinian intervals, major and minor, provided by an equitempered scale with an arbitrary number of notes. In the abscissa is represented the number of notes of the equitempered scale and in the ordinate a global parameter $\sigma$ describing the quality of the approximation to all consonant Zarlinian just intervals. $\sigma$ is calculated as the summation over all just harmonic intervals of the quadratic difference between a just harmonic interval and the nearest interval produced by the corresponding equitempered scale. The minima noted in the abscissa correspond to proposed musical scales  \cite{Cartwright2002}.
}
  \label{fig:sigma}
  \end{center}
\end{figure}

Today most musical scales are equitempered, not just, which means that the frequency interval between every pair of adjacent notes has the same ratio. The long historical battle between just and equal temperament, many of whose protagonists form part of our present story, is nonetheless one that we shall not enter here, because it takes us away from our purpose. We shall simply comment that, in mathematical terms, the problem is that there can be no one perfect musical scale because the equation $(3/2)^n= 2^m$, representing  the circle of fifths, has no nonzero integer solution.

We may consider the
numbers of notes that correspond to the optimal division of an equitempered scale into an arbitrary number of notes in such a way as best to approximate the harmonic musical intervals established by Zarlino; see Table~\ref{table:intervals} and Figure~\ref{fig:sigma}. As shown, for example, in Lattard \cite{Lattard1988} (p.\ 41) the equal division into 53 intervals given by Holder and Mercator --- anticipated in the 1st century BCE by the Chinese music theorist Ching Fang \cite{mcclain} ---  is excellent for comparing the Pythagorean and Zarlino scales; this division was also proposed and studied by Newton as being optimal among those of equal temperament. Other optimal values seen in Figure~\ref{fig:sigma} have been proposed and used: for example, that of 19 by Costeley, one of 31 intervals by Huygens, and that of 41 by Janko \cite{Lattard1988}. And there is Zarlino's 19-note harpsichord \cite{Wraight2003}.

\begin{figure}
  \begin{center}
 \includegraphics*[width= 0.6\textwidth]{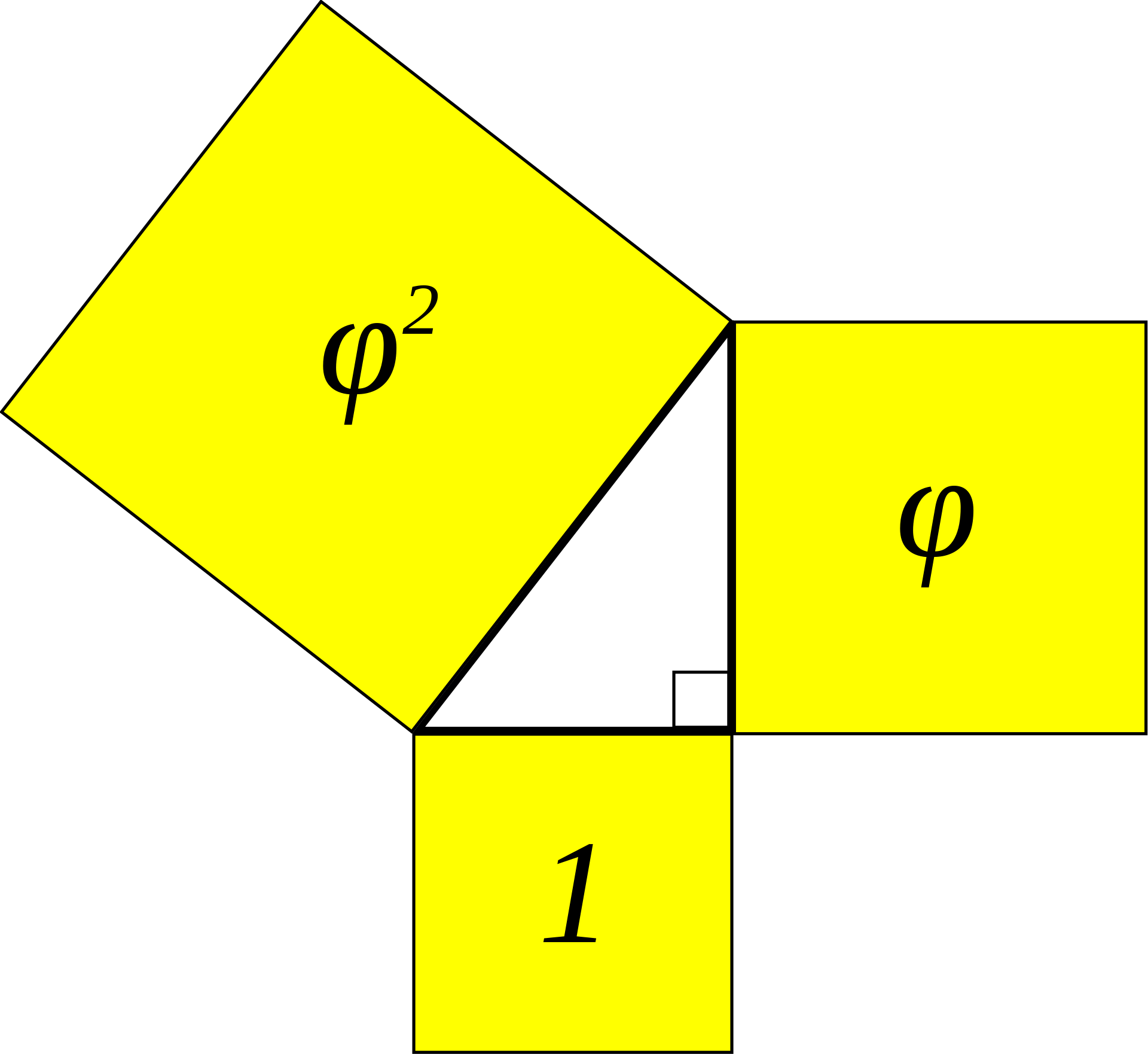}\\
  \caption{A Kepler triangle is a right triangle formed by three squares with areas in geometric progression according to the golden ratio $\varphi$, so that, by Pythagoras' theorem, $1^2+\sqrt\varphi^2=\varphi^2$ \cite{herz1993}.
}
  \label{fig:kepler_triangle}
  \end{center}
\end{figure}

If any interval is divided into arithmetic and harmonic means, the relationships formed between these means and the bounds of the closest intervals are the same:
$$
\frac{a}{\langle a,b \rangle} = \frac{|a,b|}{b}.
$$
For example, the Pythagorean intervals within the octave satisfy this property: the fifth is equivalent to that of the fourth if reversed compared to the octave, and vice versa.
The golden section, $\varphi^2 = \varphi +1$, also satisfies this property: in a series of four terms generated by the golden ratio, the two intermediate or central terms are given by the harmonic and arithmetic means of the extremes; if the geometric series of four terms are given by the proportion of the golden ratio $\varphi$:
$$
\varphi^{-1} ; 1; \varphi ; \varphi^2 ,
$$
the arithmetic and harmonic means of the extremes are
$$
\langle\varphi^{-1},\varphi^2\rangle = \frac{\varphi^{-1}+\varphi^2}{2}=\frac{\varphi^{-1}+\varphi+1}{2}= \frac{2\varphi}{2}=\varphi,
$$
$$
|\varphi^{-1},\varphi^2| = \frac{2\varphi^{-1} \varphi^2}{\varphi^{-1}+\varphi^2}=\frac{2\varphi}{2\varphi}=1,
$$
which therefore coincide with the central terms. We must return to Kepler once more to note the relation to a Kepler triangle $1^2+\sqrt\varphi^2=\varphi^2$ formed by 3 squares in geometric progression $1$; $\varphi$ ; $\varphi^2$ (Fig.~\ref{fig:kepler_triangle}): for two positive real numbers, their arithmetic mean, geometric mean, and harmonic mean are the lengths of the sides of a right triangle if and only if that triangle is a Kepler triangle \cite{herz1993}.

\begin{figure}
  \begin{center}
 \includegraphics*[width= 0.7\textwidth]{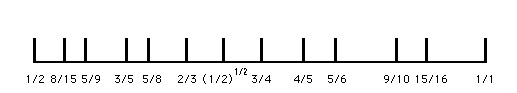}\\
  \caption{Golden scale of 12 notes; this scale is obtained on the basis of the continued fraction expansion of the golden section until the fifth convergent $8/5$, which is the quotient of two consecutive Fibonacci numbers, and of palindromic symmetry properties. Compare, for example, the solution of Newton \cite{Gouk1988}.
 The golden mean, $\varphi$, has the continued-fraction expansion 
$\displaystyle
\varphi=\frac{\sqrt{5}+1}{2}= 1+
\frac{1}{1+\displaystyle\frac{1}{1+\displaystyle\frac{1}{1+\ldots}}},
$
and the best rational approximations to $\varphi$ are given by the convergents of 
this infinite continued fraction, arrived at by cutting it off at different 
levels in the expansion: $1/1$, $2/1$, $3/2$, $5/3$, $8/5$, $13/8$, and so on; 
the convergents of the golden mean are ratios of successive Fibonacci numbers \cite{Cartwright2002}.}
  \label{fig:golden}
  \end{center}
\end{figure}

Some years ago, we took advantage of this equivalence we have highlighted above, and to the known properties of the golden section in relation to nonlinear dynamics,  to put forward a set of scales that we thought might be musically interesting, the golden scales \cite{Cartwright2002}.
The construction of golden scales is based solely on the properties of the golden section, the convergent continued fraction expansion of which coincides with the ratios of successive Fibonacci numbers, and palindromic symmetry. Figure~\ref{fig:golden} shows the construction of the 12-note golden scale. The notes reproduce all the intervals of the major and minor Zarlinian scales. From a mathematical point of view one can say that the correspondence between the Zarlinian and golden scale intervals is due to the inherent properties of palindromy of divisions in the arithmetic and harmonic means. The only difference is given by the augmented fourth that in the case of golden scale corresponds to the geometric mean of the range of an octave. The difficulties of dealing with this interval in musical terms --- the difficulty in classifying it as either consonant or dissonant --- are clearly evidenced by the name assigned to it: the \emph{Diabolus in Musica} or tritone. In the golden scale of 12 notes it is the only irrational interval. This division is necessary in the construction of a chromatic scale given that otherwise there would remain too large an interval, compared to the major semitone, between the fourth and fifth. 

The importance of results that show the equivalence between the approach of  proportions and the golden section is owing to the fundamental role that the latter plays as the most irrational number  in the description of the behaviour of nonlinear dynamical systems \cite{percival1987}. There is a hierarchy among irrational numbers according to how difficult they are to approximate with rationals; it is in this sense that one irrational is more irrational than another, and this difficulty of approximation, which is related to the continued fraction expansion mentioned in Figure~\ref{fig:golden}, is given by Hurwitz's theorem \cite{hurwitz1891}.

\section*{Nonlinear dynamics, resonance, and synchronization}

\begin{figure}
  \begin{center}
 \includegraphics*[height= 0.4\textwidth]{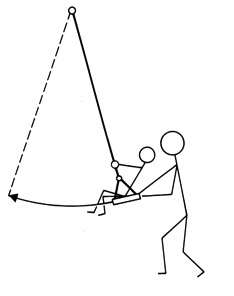} 
 \includegraphics*[height= 0.4\textwidth]{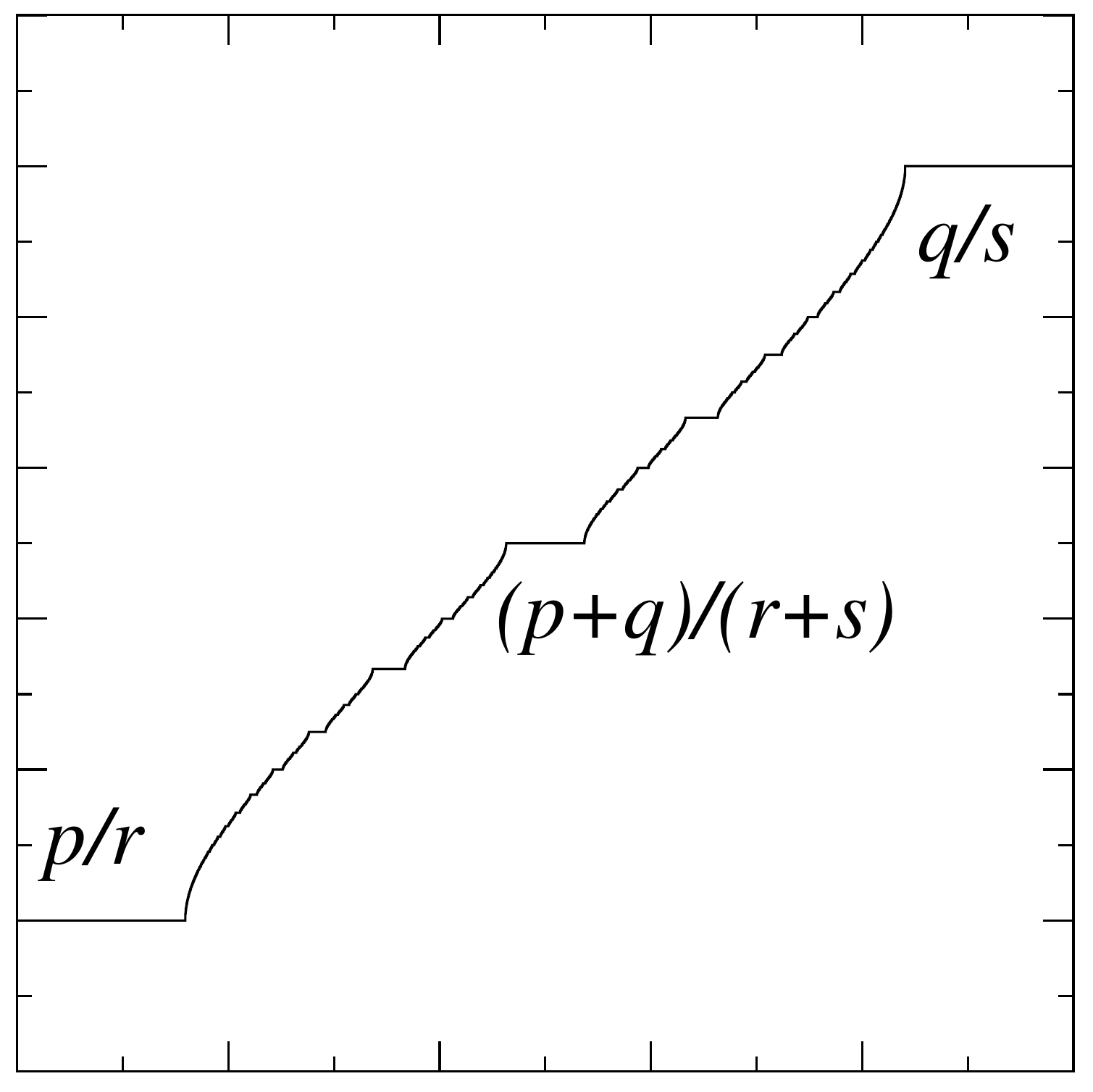} \\
 \includegraphics*[height= 0.2\textwidth]{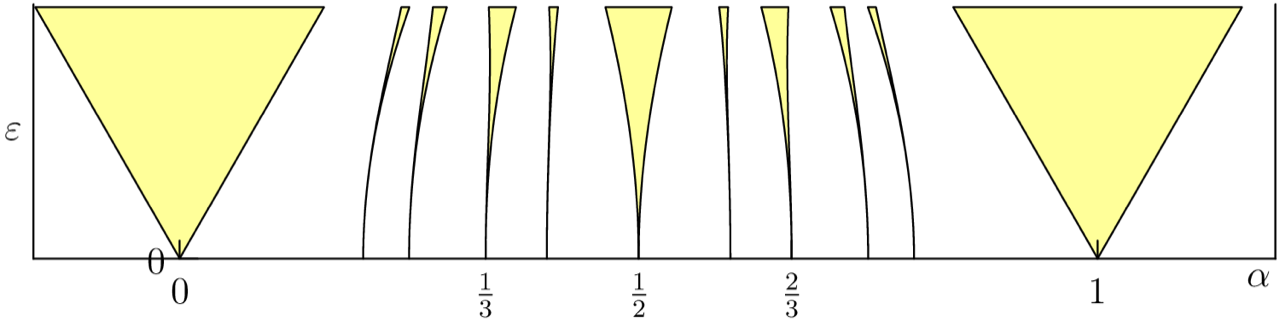} 
  \caption{a) The simplest dynamical system: the simple pendulum, here represented by a swing forced with periodic impulses.
b) Dynamics of a periodically forced swing: the regions of synchronization, known as Arnol'd tongues. Each region is born in a rational number on the horizontal axis where  the period of the impulses that act on the swing is represented. The denominator is the number of complete oscillations of the swing that are performed during a number of impulses coincident with the numerator of this rational number. The vertical axis represents the intensity of impulses applied.c) A cut along a horizontal line produces the so-called devil's staircase.
}
  \label{fig:arnold_tongues}\label{fig:swing}
  \end{center}
\end{figure}

In order to explore the connection of the golden scales with nonlinear dynamics, 
we introduce as an example one of the simplest dynamical systems: the forced pendulum, which is equivalent to a child's swing (Figure~\ref{fig:swing}).
A phenomenon familiar to everyone who has ever played on a swing  is that of synchronization: for each complete oscillation of the swing it receives a single impulse, from our legs, or from the arms of a pusher. Since we have one oscillation for each impulse, the frequencies of the swing and the impulses are equal, i.e., the ratio of frequencies is $1/1 = 1$. The phenomenon of synchronization, first discovered and analysed by Huygens in pendulum clocks in 1665 \cite{Huygens1665}, is one of the most universal examples of dynamical behaviour and can be observed at all levels of the physical world.
For example, in celestial mechanics, we observe always the same face of the Moon from the Earth, as it takes the Moon exactly the same time to turn on itself as to make a revolution around the Earth. In other words, the Moon's rotational and orbital periods are synchronized.

\begin{figure}
  \begin{center}
 \includegraphics*[width= 0.6\textwidth]{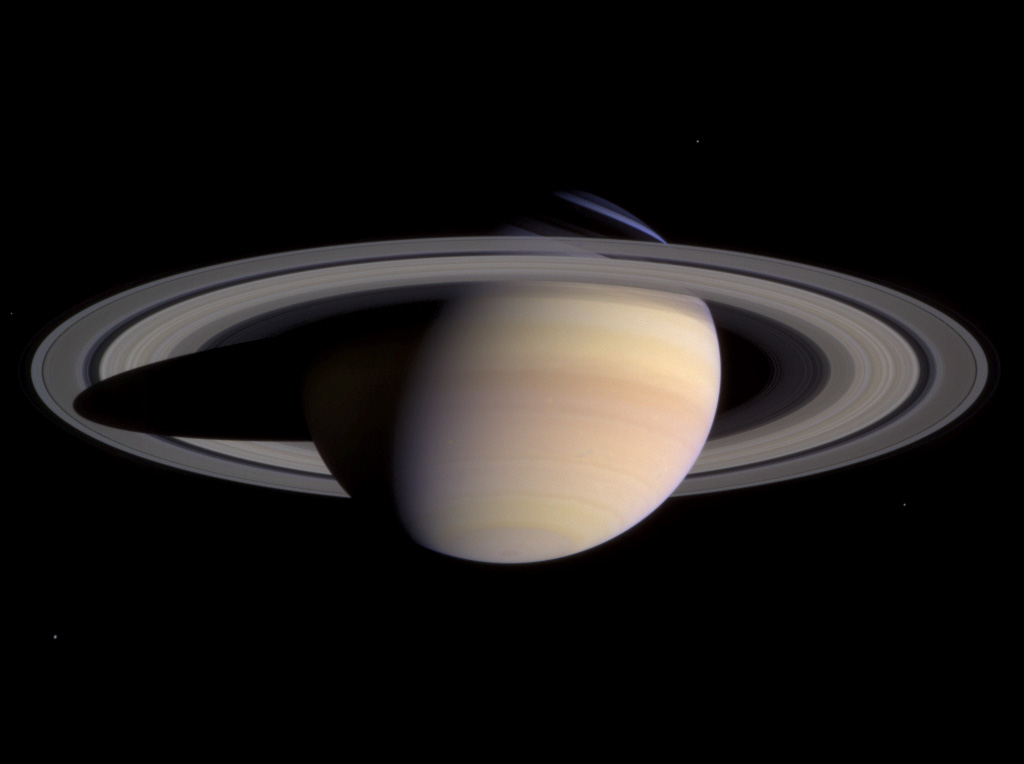} \\
  \caption{The rings of Saturn, made up of a myriad of highly-reflective blocks of ice, have resonances (ring gaps), similar to musical intervals, within which material is almost absent. The most prominent is the Cassini division, the thick dark band in this image.
}
  \label{fig:saturn}
  \end{center}
\end{figure}

Returning to the swing, in Figure~\ref{fig:arnold_tongues} we have a graph that represents all the dynamical behaviour of a swing   possible when varying the intensity of the impulses (vertical axis) and their frequency (horizontal axis). There are regions of synchronization more complicated than $1/1$, for example, $1/2$ or $2/3$, that have the shape of a tongue and which are, in fact, called Arnol'd tongues in honour of the Russian mathematician who studied them \cite{Arnold1978}. We may note that at every rational number an Arnol'd tongue is born. 
Such more complex synchronizations can be seen across the natural world, including in celestial mechanics \cite{Escribano2008,Vanyo2011}.  Examples in our Solar System are Mercury, tidally locked with the Sun in a 3/2 spin--orbit resonance --- a species of celestial fifth --- the absence of asteroids in the Kirkwood gaps in the asteroid belt owing to resonances with Jupiter, and gaps in Saturn's rings such as the Cassini division (Figure~\ref{fig:saturn}), caused by resonances with moons.

From a Pythagorean perspective Figure~\ref{fig:arnold_tongues} is very interesting. Here is a simple physical system that can distinguish rational numbers, a seemingly purely mathematical concept, from all others. If one cuts the figure by a horizontal line and represents the intervals in which a given region of synchronization is stable, we get the so-called devil's staircase \cite{Cartwright2010}. This staircase is made up of infinite steps: between two successive steps, there is always another. In addition, the staircase is self-similar, and fractal: if one zooms in to any piece of this staircase one sees that the part is equal to the full staircase (Figure~\ref{fig:arnold_tongues}c). 

The sizes of the various entrainment regions are ordered in a way related to a concept from number theory: Farey sequences. An $n$-Farey sequence is the increasing succession of rational numbers whose denominators are less than or equal to $n$.  We call two rational numbers `adjacent' if they are consecutive in the Farey sequence. A necessary and sufficient condition for $p/q$ and $r/s$ to be adjacent is $|ps  - qr |=1$. A rational number $\alpha$ belonging to the open interval $(p/q;r/s)$ where $p /q$ and $r/s$ are adjacent will be called `mediant' if there is no other rational in the interval having smaller denominator. It is known that $\alpha=(p+r)/(q+s)$ and is unique. Observation of Fig.~\ref{fig:arnold_tongues} allows us to guess that the synchronization zone characterized by a mediant number of two adjacent rationals is the greatest of all the zones situated between those rationals. In addition it has, obviously, the least period. This property is generic and not specific of the case illustrated in Fig.~\ref{fig:arnold_tongues}. As a consequence, we have that between two solutions corresponding to successive convergents of the golden section, say $n$ and $n+1$ (which are represented by Fibonacci quotients and are also Farey adjacents), the widest region corresponds to the successive convergent of the golden section, say $n + 2$, that coincides also with their Farey sum. This creates a hierarchy in parameter space that follows a Farey tree. This means also that, starting with the first two convergents, $1/1$ and $1/2$, we can obtain all the convergents of the golden section by the Farey sum operation; any convergent is the widest synchronization region between the two parent regions. Wider stability intervals means also that such solutions are more robust to parameter perturbations and, thus, that are more relevant for the generic modelling of the dynamical system under study.

\begin{figure}
  \begin{center}
 \includegraphics*[width= 0.6\textwidth]{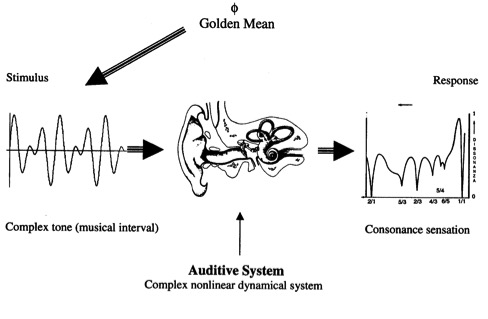} \\
  \caption{Considering our auditory system as a general nonlinear system, and observing the analogy between the phase diagram (Figure~\ref{fig:arnold_tongues}) and the consonance diagram of Plomp \cite{Plomp1965}, one might expect that the regions of synchronization can describe the phenomenon of musical consonance. However, a fundamental problem remains: the stability of the regions of synchronization.}
  \label{fig:simplified}
  \end{center}
\end{figure}

As can be seen in the diagram of  Figure~\ref{fig:simplified}, it might appear that synchronization should explain the harmonic intervals, given also the similarity of the consonance curves with those obtained, for example, through the psychoacoustic theory of Plomp \cite{Plomp1965}, which basically affirms that musical intervals seem consonant if the frequency differences between components exceed a critical bandwidth.
However, a fundamental problem presents itself: that of stability. A harmonic interval implies the presence of at least two independent frequencies (those which define the interval) and synchronization (such as that of the Earth--Moon system) becomes unstable in the presence of two incommensurable frequencies.
This is a problem faced by all theories of music using the Pythagorean or Zarlinian just intervals. If an interval is represented by a rational ratio, like the fifth, a small detuning in one of the frequencies automatically
breaks this rational ratio. Paradoxically, the smaller the detuning, the greater the whole numbers needed to describe the new rational ratio. It is clear that this mathematical instability does not correspond with our perception: the
feeling of consonance is a maximum for the ratios defined by small numbers, for example $3/2$ for the fifth, and degrades proportionally when the range is out of tune to a certain extent. 
This represents an essential drawback to fitting the existence of harmonic intervals and musical scales into a dynamical-systems paradigm.

\section*{Quasiperiodic forcing and three-frequency resonance}

\begin{figure}
  \begin{center}
\includegraphics*[height= 0.4\textwidth]{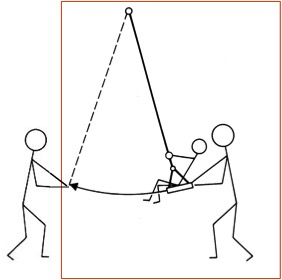}
   \includegraphics*[width= 0.8\textwidth]{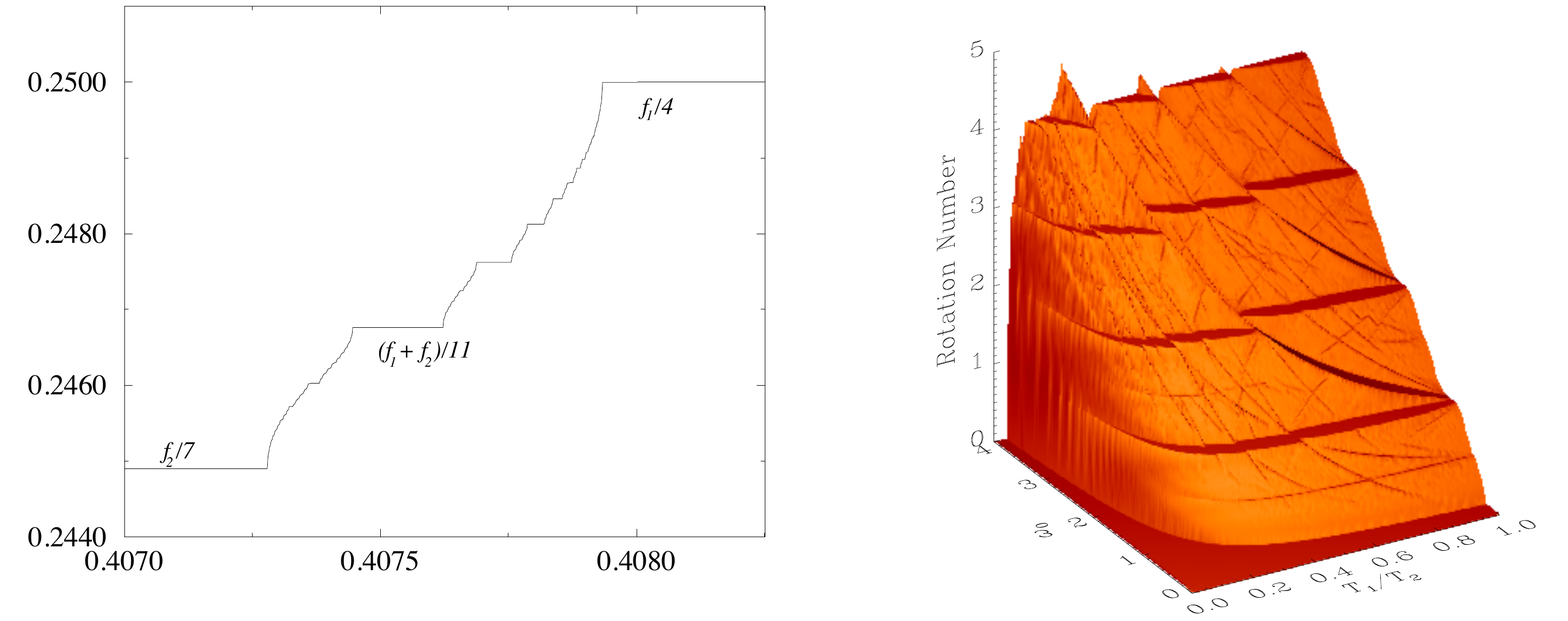}\\
    \caption{Adding a second pusher to Figure~\ref{fig:swing} gives a quasiperiodically forced swing: the presence of two simultaneous periodic forces pushing with different frequencies destroys the stability of regions of synchronization.    
With a real swing, where one generally pushes once per period of the swing, this setup would not be so easy to carry out for an arbitrary forcing frequency.
Dynamics of a quasiperiodically forced swing. 
  (a) Three-frequency devil's staircase for forcing frequencies $f_1 = 1$, $f_2 = 12/7$.  The generalized regions of synchronization represent resonances at three frequencies.
  (b) Devil's ramps: the global organization of three-frequency resonances as a function of the external frequency ratio and the intrinsic frequency \cite{Cartwright1999b}.}
   \label{fig:devil's_ramps}\label{fig:quasi-periodic}
  \end{center}
\end{figure}

Ultimately the very existence of harmonic musical intervals in the brain may be due to the properties of our auditory system, which is a highly nonlinear dynamical system. A harmonic interval is a particular case of a chord, indeed a chord with only two notes. Although a harmonic interval is an example of a very simple chord, the above-described stability problem remains. We need at least two independent parameters for describing a harmonic interval, for example, the two fundamental frequencies of the notes composing the interval. Thus, as we have sketched in Fig.~\ref{fig:simplified}, we can consider our auditory system as a dynamical system in this case forced with two independent frequencies. 
Consequently, in order to tackle the problem of harmonic intervals one needs to consider a slightly more complex dynamical system than the periodically forced swing: a quasiperiodically forced swing. The system and its responses are  represented in Figure~\ref{fig:devil's_ramps}. In this case we have regions of generalized synchronization. Within these regions a relation between three resonant frequencies is satisfied.
In the case of simple synchronization there is a relationship between two resonant frequencies that can be described mathematically as
$$
p f_1 + q f_2 =0 ,
$$
where $p$, $q$ are integers, $f_1$ the forcing frequency and $f_2$ the frequency of the response.
It is easy to generalize this condition to three frequencies
$$
p f_1 + q f_2 + r f_3 =0 ,
$$
where $p$, $q$, $r$, are integers, $f_1$, $f_2$, the forcing frequencies, and $f_3$, the frequency of the response.
As is the case for simple synchronization, a section of the phase portrait, shown in Figure~\ref{fig:devil's_ramps}a, shows the existence of a three-frequency devil's staircase shown in Figure~\ref{fig:devil's_ramps}b. The aspect that it is important to note about these three-frequency resonances is that they are also organized by means of hierarchical rules dictated by number theory, i.e., in a Pythagorean fashion \cite{Cartwright1999b}.

\section*{Dynamical systems in man and the heavens: The examples of the auditory system (the microcosmos) and celestial mechanics (the macrocosmos)}

\begin{figure}
  \begin{center}
 \includegraphics*[width= 0.6\textwidth]{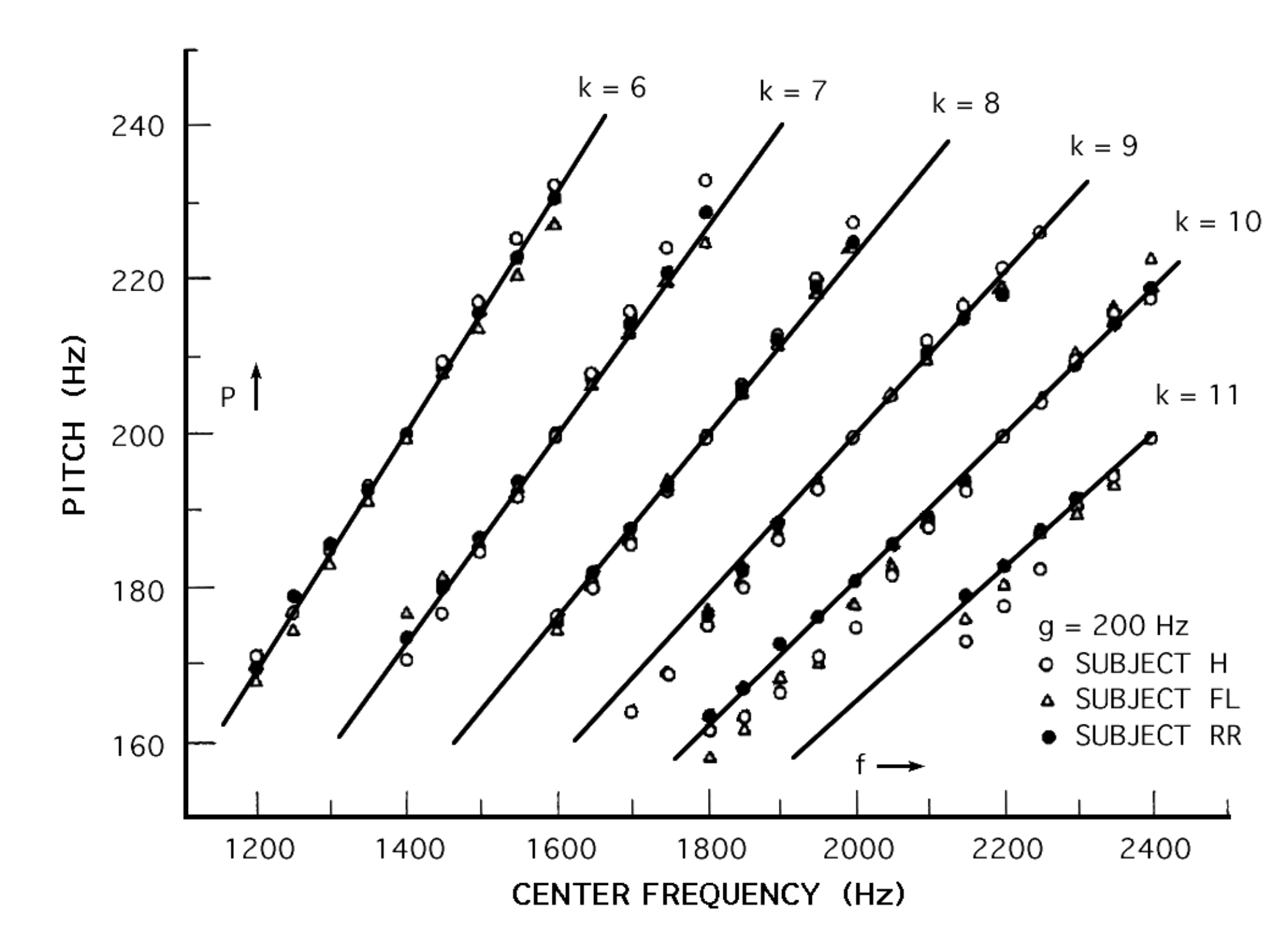}\\
  \caption{The solid lines represent the theoretical solutions provided by three-frequency resonances to the key psychophysical phenomenon of pitch shift of the missing fundamental presented in Figure~\ref{fig:residue}. The various symbols represent the results of psychoacoustic experiments with three different human different subjects. As can be seen in the figure, the agreement between theory and experiment is rather good \cite{Cartwright1999a}.}
  \label{fig:pitch_shift}
  \end{center}
\end{figure}

Resonances of three frequencies are not simply a mathematical curiosity but rather are a feature found at different levels of natural phenomena. 
In the Solar System, the behaviour of some bodies in the Kuiper belt and in Saturn's rings is described by three-frequency resonances \cite{Murray2001}. 
In addition, for music the organization of these resonances corresponds with the perception of the missing fundamental and accurately describes the phenomenon of pitch shift, as can be seen in Figure~\ref{fig:pitch_shift} \cite{Cartwright1999a,Cartwright2001}.

\begin{figure}
  \begin{center}
 \includegraphics*[width= 0.6\textwidth]{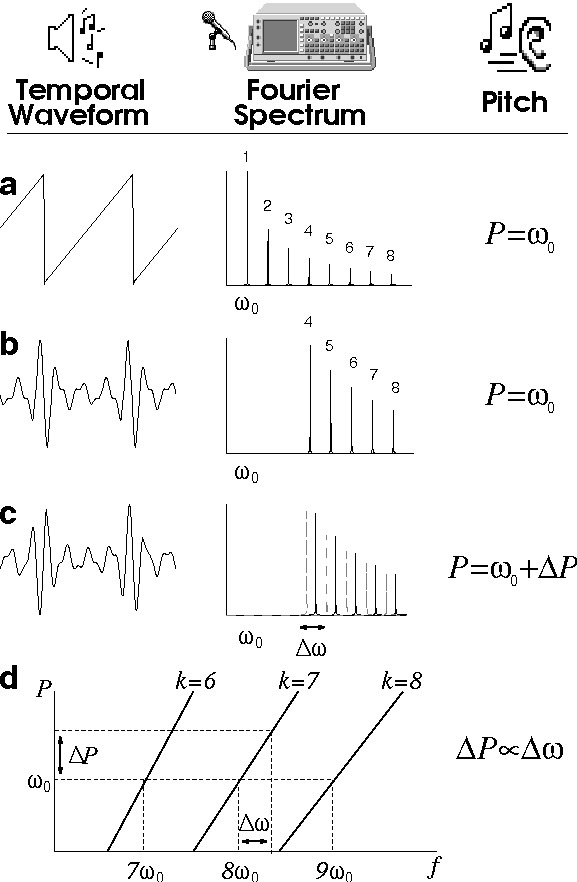}\\
  \caption{
  The phenomenon of the missing fundamental or residue:
  Fourier spectra and pitches of complex tones. Whereas pure tones have a
sinusoidal waveform corresponding to a single frequency, almost all musical 
sounds are complex tones that consist of a lowest frequency component, or 
fundamental, together with higher frequency overtones. The fundamental plus 
overtones are together collectively called partials.
{a}) A harmonic complex tone. The overtones are successive integer 
multiples $k=2,3,4\ldots$ of the fundamental $\omega_0$ that determines the 
pitch. The partials of a harmonic complex tone are termed harmonics.
{b}) Another harmonic complex tone. The fundamental and the first few higher
harmonics have been removed. The pitch remains the same and equal to the
missing fundamental. This pitch is known as virtual or residue pitch.
{c}) An anharmonic complex tone. The partials, which are no longer harmonics,
are obtained by a uniform shift $\Delta\omega$ of the previous harmonic case 
(shown dashed). Although the difference combination tones between successive 
partials,
$\omega_C=\omega_2-\omega_1$, remain unchanged and equal to the missing fundamental, the pitch 
shifts by a quantity $\Delta P$ that depends linearly on $\Delta\omega$.
{d}) Pitch shift. Pitch as a function of the central frequency 
$f=(k+1)\omega_0+\Delta\omega$ 
of a three-component complex tone $\{k\omega_0+\Delta\omega, 
(k+1)\omega_0+\Delta\omega, (k+2)\omega_0+\Delta\omega\}$. The pitch-shift 
effect is shown here for $k=6$, $7$, and $8$. Three-component complex
tones are often used in pitch experiments because they elicit a clear residue 
sensation and can easily be obtained by amplitude modulation of a pure tone of 
frequency $f$ with another pure tone of frequency $\omega_0$. When $\omega_0$ 
and $f$ are rationally related, $\Delta\omega=0$, and the three frequencies are 
successive multiples of some missing fundamental. At this point $\Delta P=0$,
and the pitch is $\omega_0$, coincident with the frequency of the missing 
fundamental \cite{Cartwright1999a}.
}
  \label{fig:residue}
  \end{center}
\end{figure}

Suppose that a periodic signal be presented to the ear. The pitch of the signal
can be quantitatively well described by the frequency of the fundamental, say
$\omega_0$; see Fig.~\ref{fig:residue}(a). The number of harmonics and their
relative amplitudes gives the timbral characteristics to the sound. Now suppose that the fundamental and  some
of the first few higher harmonics are removed,  Fig.~\ref{fig:residue}(b).
Although the timbral sensation changes, the pitch of the complex remains
unchanged and equal to the missing fundamental. This is termed residue perception. Moreover, psychophysical experiments can be done by shifting the remaining partials, as shown in Fig.~\ref{fig:residue}(c), whereon the perceived pitch also shifts  in accordance with three-frequency resonances, as sketched in Fig.~\ref{fig:residue}(d), and shown in psychoacoustic experiments in Figure~\ref{fig:pitch_shift}.
Thus, the perception of the pitch of complex sounds, and in particular the key phenomenon of perception of the missing fundamental, or residue perception, can be described by means of the dynamical properties of our auditory system \cite{Cartwright1999a,Cartwright2001}. This represents a key step in the development of a nonlinear approach for musical perception.

It has also been suggested that the residue is one of the fundamental phenomena of musical perception and the basis of the consonant structure of musical intervals \cite{Terhardt1974,Roederer1995}. Historically, the theories of Rameau of the fundamental sound \cite{Rameau1772} and that of Tartini based on the third sound \cite{Tartini1754} represent anticipations of this perceptual and nonlinear origin of musical harmony.
In the case of chords, the residue can perhaps be identified with the fundamental sound of Rameau ``always the lowest and deepest part'' and, therefore, becomes the basis of the development of harmonious music. As Rameau put it
\begin{quote}
the harmony of these consonances can be perfect only if the first sound is found below them, serving as their base and fundamental. \cite{Rameau1772}
\end{quote} 
Moreover, we may quote from the composer and violinist Tartini's description of the third sound:
\begin{quote}
Let the following intervals be played perfectly by a violinist simultaneously with a strong, sustained bow. A third sound will be heard ...  The same will happen if the presented intervals are played by two violin players five or six steps apart, each playing his note at the same time, and always with a strong, sustained bow. A listener in the middle of the two players will hear this third sound much more ... \cite{Tartini1754}
\end{quote} 
We think that Tartini is describing  residue perception and not simply a combination tone here. In the second experimental situation described, the intensity in one ear due to the opposite violin should be --- owing to the auditory shadow of the head --- very low, and thus not able to produce a detectable combination tone. It is very probable that he is describing the perception of the residue in a dichotic situation; that of a different stimulus in each ear. 
In this way, the rules of musical composition advanced by Rameau and Tartini find a physical rationale in the theory of dynamical systems.

These results suggest that there is something fundamental in harmonic musical intervals and bring us back to the original musical problem treated by the Pythagoreans: what is the origin of these numerical relations that describe so precisely the basic elements of aesthetics of auditory perception? Historical and archaeological research likewise leads one to consider that the existence of privileged intervals may be inherent to the physiology of the hearing system rather than being due to culture. For example, analysis of a bone flute found in a Neanderthal site in Slovenia \cite{Anon1997,Fink2003} and other, still functional, prehistoric flutes from China \cite{Zhang1999} indicate that these instruments already produced consonant musical intervals. These findings are notable because not only do they highlight that musical behaviour may not be exclusive to \emph{Homo sapiens}, and thus to human culture, but also that harmonic musical intervals are a phenomenon existing at the base of music and therefore may reflect the existence of universal physiological mechanisms in the auditory system of mammals and even perhaps of all life.

Moreover, such a foundation prompts a revisiting of the Pythagorean conception of music as a paradigm of
cosmological order and the role of numbers in such a pattern. In this paradigm, there are of course many other possible dynamical responses; for example, very complex or erratic behaviour, similar to what is observed in random or stochastic processes. That phenomenon is known as deterministic chaos and was first identified by the meteorologist Edward Lorenz in a highly simplified climate model \cite{Lorenz1963} that he analysed with one of the first electronic computers. An example of this behaviour in the Solar System is provided by the moon of Saturn, Hyperion, which rotates in just such a chaotic way \cite{Wisdom1984}. Of course this deterministic chaos is not the same thing as Pythagorean chaos; yet the modern term is a deliberate echo of the ancient Greeks \cite{yorke1975}.
The widespread presence of synchronization and generalized synchronization, as well as chaos, in planetary motions vindicates a key Pythagorean explanation of cosmological order.

\section*{Musica Universalis?}

\begin{figure}
  \begin{center}
 \includegraphics*[width= 0.6\textwidth]{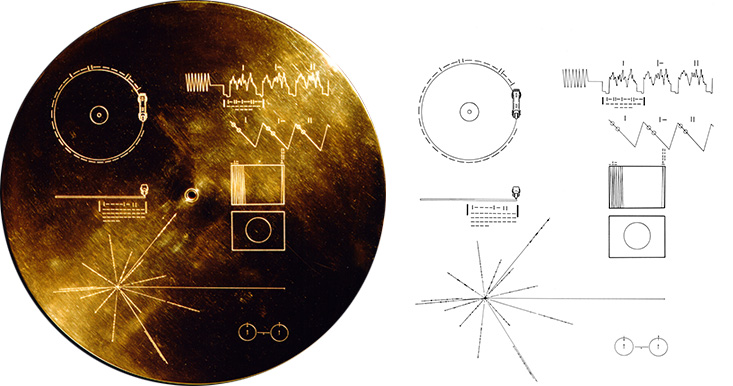}\\
  \caption{All is number, \emph{numeri regunt mundum}, was the Pythagorean ethos. The cover of the Voyager golden records encode the information about where Earth is in the universe, and about how to play the record and extract the information from it, using numbers 
  \cite{sagan1978,pescovitz2017}.
  Image: NASA/JPL.}
  \label{fig:golden_record}
  \end{center}
\end{figure}

The only manmade object so far to leave our Solar System and enter interstellar space  is Voyager 1, launched in 1977 and now navigating beyond the heliosphere, over 140 times further from the Sun than we are on Earth.
Aboard both Voyager 1 and its twin, Voyager 2, which is nearly as far away from us, but headed in a completely different direction in space, there is music. The fortuitously but aptly named Golden Record (Fig.~\ref{fig:golden_record}) --- the designation comes from the record's material composition; it is gold plated copper --- carried by the two Voyager spacecraft was put together by a group of people led by Carl Sagan  to be sent into space like a message in a bottle. As Sagan wrote \cite{sagan1978}
\begin{quote}
I was delighted with the suggestion of sending a record ... we could send music. ...  Perhaps a sufficiently advanced civilization would have made an inventory of the music of species on many planets and, by comparing our music with such a library, might be able to deduce a great deal about us. ... Because of the relation between music and mathematics, and the anticipated universality of mathematics, it may be that much more than our emotions are conveyed by the musical offering on the Voyager record.
\end{quote}
Sagan took some of these ideas from work of von Hoerner \cite{hoerner1974}, who concluded that
\begin{quote}
Some other civilisations in space may have no music at all, for various biological or mental reasons. Some others may have vastly different things which they call music but which are incomprehensible to us, for similar reasons. But it seems that some of our basic musical principles are universal enough to be expected at a good fraction of other places, too: a chromatic scale of exactly 5, 12, or 31 equal parts, from which rather arbitrary numbers and sequences can be selected for melodic scales.
\end{quote}
It has long been an idea that music may be a universal language, and the logical extension to that idea is that music is likely to be a common channel of communication with alien intelligences.  Huygens \cite{huygens1698}  wrote on music and musical scales in extraterrestrial civilizations in his book \emph{Cosmotheoros}, from which we think it worthwhile to quote his ideas at some length:
\begin{quote}
It's the same with Musick as with Geometry, it's every where immutably the same, and always will be so. For all Harmony consists in Concord, and Concord is all the World over fixt according to the same invariable measure and proportion. So that in all Nations the difference and distance of Notes is the same, whether they be in a continued gradual progression, or the voice makes skips over one to the next. Nay very credible Authors report, that there's a sort of Bird in America, that can plainly sing in order six musical Notes: whence it follows that the Laws of Musick are unchangeably fix'd by Nature, and therefore the same Reason holds valid for their Musick, as we e'en now proposed for their Geometry. For why, supposing other Nations and Creatures, endued with Reason and Sense as well as we, should not they reap the Pleasures arising from these Senses as well as we too? I don't know what effect this Argument, from the immutable nature of these Arts, may have upon the Minds of others; I think it no inconsiderable or contemptible one, but of as great Strength as that which I made use of above to prove that the Planetarians had the sense of Seeing.

But if they take delight in Harmony, 'tis twenty to one but that they have invented musical Instruments. For, if nothing else, they could scarce help lighting upon some or other by chance; the sound of a tight String, the noise of the Winds, or the whistling of Reeds, might have given them the hint. From these small beginnings they perhaps, as well as we, have advanced by degrees to the use of the Lute, Harp, Flute, and many string'd Instruments. But altho the Tones are certain and determinate, yet we find among different Nations a quite different manner and rule for Singing; as formerly among the Dorians, Phrygians, and Lydians, and in our time among the French, Italians, and Persians. In like manner it may so happen, that the Musick of the Inhabitants of the Planets may widely differ from all these, and yet be very good. But why we should look upon their Musick to be worse than ours, there's no reason can be given; neither can we well presume that they want the use of half-notes and quarter-notes, seeing the invention of half-notes is so obvious, and the use of 'em so agreeable to nature. Nay, to go a step farther, what if they should excel us in the Theory and practick part of Musick, and outdo us in Consorts of vocal and instrumental Musick, so artificially compos'd, that they shew their Skill by the mixtures of Discords and Concords? and of this last sort 'tis very likely the 5th and 3d in use with them.

This is a very bold Assertion, but it may be true for ought we know, and the Inhabitants of the Planets may possibly have a greater insight into the Theory of Musick than has yet bin discover'd amongst us. 
\end{quote}
For just such reasons Sagan and his collaborators sent music to whoever in the future might recover the wandering spacecraft.  They sent Bach, Beethoven, and Mozart, but also music from around the world \cite{pescovitz2017}. One piece that they included, prefacing a selection of the natural and human sounds of Earth and commissioned just for the Golden Record, explicitly encodes Pythagorean ideas:
\begin{quote}
the giddy whirl of tones reflecting the motions of the Sun's planets in their orbits --- a musical readout of Johannes Kepler's \emph{Harmonica Mundi}, the sixteenth-century mathematical tract whose echoes may still be found in the formulas that make Voyager possible. Kepler's concept was realized on a computer at Bell Telephone Laboratories by composer Laurie Spiegel in collaboration with Yale professors John Rogers and Willie Ruff. Each frequency represents a planet; the highest pitch represents the motion of Mercury around the Sun as seen from Earth; the lowest frequency represents Jupiter's orbital motion. Inner planets circle the Sun more swiftly than the outer planets. The particular segment that appears on the record corresponds to very roughly a century of planetary motion. Kepler was enamored of a literal ``music of the spheres,'' and I think he would have loved their haunting representation here.
\cite{sagan1978}
\end{quote}
We have argued here that Pythagorean ideas of a Musica Universalis are correct owing to the universality of resonances of nonlinear dynamical systems in the cosmos. If that conception is true then this Keplerian piece may act for its discoverers as a Rosetta Stone that unlocks for them the rest of Earth music.

\section*{Coda}

We have presented here our vision of how current-day research in dynamical systems is a continuation of a millenary tradition  linking mathematics, celestial mechanics, and music stretching back to Pythagoras. We continue to work on aspects of dynamical systems applied to music perception.
The nervous system is of great complexity, both functional and structural, and in particular our auditory system is wired in an extraordinarily precise and complex way. We are currently working on models that take into account the physiology of the auditory system in order to explain some aspects of musical perception that may reflect the existence of an underlying dynamics.
After applying the dynamics of three-frequency resonances to musical pitch and the theory of residue perception, we now seek a nonlinear theory for musical consonance based on the fundamental role of nonlinear dynamics in determining the physiological and psychoacoustic responses of the auditory system to stimulation by complex sounds.

\bibliographystyle{spmpsci}      
\bibliography{pythagoras_revisited}

\end{document}